\documentclass[12pt,a4paper]{amsart}
\usepackage{amssymb}

\def\bfit{\bfseries\itshape}

\input xypic
\xyoption{all}
\xyoption{arc}

\def\equat{\refstepcounter{theo}$$~}
\def\endequat{\leqno{\boldsymbol{(\arabic{theo})}}~$$}

\def\Gb{{\mathbf G}} 

\def\Pb{{\mathbf P}}
\def\Tb{{\mathbf T}}

\def\Bb{{\mathbf B}}
\def\Xb{{\mathbf X}}

\def\Xbh{{\hat{\Xb}}}
\def\PCB{{\boldsymbol{\mathcal{P}}}}
\def\BCB{{\boldsymbol{\mathcal{B}}}}
\def\OCB{{\boldsymbol{\mathcal{O}}}}

\def\d{\delta}

\def\t{\tau}

\def\to{\rightarrow}
\def\longto{\longrightarrow}

\def\vide{\varnothing}
\def\fin{~${\scriptstyle{\blacksquare}}$}

\def\lexp#1#2{\kern\scriptspace\vphantom{#2}^{#1}\kern-\scriptspace#2}
\def\le{\hspace{0.1em}\mathop{\leqslant}\nolimits\hspace{0.1em}}
\def\ge{\hspace{0.1em}\mathop{\geqslant}\nolimits\hspace{0.1em}}

\begin{document}

\title{On the irreducibility of Deligne-Lusztig varieties}

\author{C\'edric Bonnaf\'e \& Rapha\"el Rouquier}
\address{C\'edric Bonnaf\'e: 
Laboratoire de Math\'ematiques de Besan\c{c}on (CNRS: UMR 6623), 
Universit\'e de Franche-Comt\'e, 16 Route de Gray, 25030 Besan\c{c}on
Cedex, France} 

\makeatletter
\email{bonnafe@math.univ-fcomte.fr}
\urladdr{http://www-math.univ-fcomte.fr/pp\_Annu/CBONNAFE/}

\address{Rapha\"el Rouquier: 
Department of Pure Mathematics, University of Leeds, Leeds LS2 9JT, U.K.}
\email{rouquier@maths.leeds.ac.uk}
\urladdr{http://www.maths.leeds.ac.uk/{$\sim$}rouquier/}

\makeatother


\date{\today}

\begin{abstract} 
Let $\Gb$ be a connected reductive algebraic
group defined over an algebraic closure 
of a finite field and let $F : \Gb \to \Gb$ be an endomorphism such that $F^\d$ 
is a Frobenius endomorphism for some $\d \ge 1$. Let $\Pb$ be a parabolic
subgroup of $\Gb$. We prove that the Deligne-Lusztig 
variety $\{g\Pb~|~g^{-1}F(g)\in \Pb\cdot F(\Pb)\}$ is irreducible 
if and only if $\Pb$ is not contained in a proper $F$-stable parabolic 
subgroup of $\Gb$.
\end{abstract}

\maketitle

\pagestyle{myheadings}

\markboth{\sc C. Bonnaf\'e \& R. Rouquier}{\sc Irreducibility of 
Deligne-Lusztig varieties}

Let $\Gb$ be a connected reductive group over an algebraic closure of
a finite field and let 
$F : \Gb \to \Gb$ be an endomorphism such that some power of $F$ is a Frobenius 
endomorphism of $\Gb$. If $\Pb$ is a parabolic subgroup of $\Gb$, we set
$$\Xb_\Pb=\{g\Pb \in \Gb/\Pb~|~g^{-1} F(g) \in \Pb\cdot F(\Pb)\}.$$
This is the Deligne-Lusztig variety associated to $\Pb$. The aim 
of this note is to prove the following result:

\bigskip

\noindent{\bf Theorem A.} 
{\it Let $\Pb$ be a parabolic subgroup of $\Gb$. Then 
$\Xb_\Pb$ is irreducible if and only if $\Pb$ is not contained 
in a proper $F$-stable parabolic subgroup of $\Gb$.}

\bigskip

Note that this result has been obtained independently by
Lusztig (unpublished) and Digne and Michel 
\cite[Proposition 8.4]{dimi} in the case 
where $\Pb$ is a Borel subgroup: both proofs are obtained by counting 
rational points of $\Xb_\Pb$ in terms of the Hecke algebra. We present here 
a geometric proof (inspired by an argument of Deligne 
\cite[proof of Proposition 4.8]{lucox}) which reduces the problem to
the irreducibility of the Deligne-Lusztig variety associated to 
a Coxeter element: this case has been treated by Deligne and Lusztig 
\cite[Proposition 4.8]{lucox}. 

\smallskip
Before starting the proof of this Theorem, we first describe an equivalent 
statement. Let $\Bb$ be an $F$-stable Borel subgroup of $\Gb$, let $\Tb$ be an $F$-stable 
maximal torus of $\Bb$, let $W$ be the Weyl group of $\Gb$ relative to $\Tb$
and let $S$ be the set of simple reflections of $W$ with respect to $\Bb$. We 
denote again by $F$ the automorphism of $W$ induced by $F$. 
Given $I \subset S$, let $W_I$ denote the standard parabolic subgroup of $W$ 
generated by $I$ and let $\Pb_I=\Bb W_I \Bb$. 
We denote by $\PCB_I$ the variety of parabolic subgroups of $\Gb$ 
of type $I$ (i.e. conjugate to $\Pb_I$) and by $\BCB$ the variety of Borel 
subgroups of $\Gb$ (i.e. $\BCB=\PCB_\vide$). For $w \in W$, we denote by 
$\OCB_I(w)$ the $\Gb$-orbit of $(\Pb_I,\lexp{w}{\Pb_{F(I)}})$ in $\PCB_I \times \PCB_{F(I)}$. 
Note that $\OCB_I(w)$ depends only on the double coset $W_I w W_{F(I)}$. We
define now
$$\Xb_I(w)=\{\Pb \in \PCB_I~|~(\Pb,F(\Pb)) \in \OCB_I(w)\}.$$
The group $\Gb^F$ acts on $\Xb_I(w)$ by conjugation. We set 
$\OCB(w)=\OCB_\vide(w)$ and $\Xb(w)=\Xb_\vide(w)$.

\bigskip

\noindent{\bf Theorem A'.} 
{\it Let $I \subset S$ and let $w \in W$. 
Then $\Xb_I(w)$ is irreducible if and only if $W_Iw$ is not contained in 
a proper $F$-stable standard parabolic subgroup of $W$.}

\bigskip
\def\Qb{{\mathbf{Q}}}

\noindent{\sc Remark 1 - } 
Let us explain why the Theorems A and A' are equivalent. Let $\Pb_0$ be a parabolic 
subgroup of $\Gb$. Let $I$ be its type and let $g_0 \in \Gb$ be such that 
$\Pb_0=\lexp{g_0}{\Pb_I}$. Let $w \in W$ be such that 
$g_0^{-1} F(g_0) \in \Pb_I w \Pb_{F(I)}$. The pair $(I,W_I w W_{F(I)})$ 
is uniquely determined by $\Pb_0$. Then, the map 
$\Xb_{\Pb_0} \to \Xb_I(w)$, $g \Pb_0 \mapsto \lexp{gg_0}{\Pb_I}$ is an isomorphism 
of varieties (indeed, it is straightforward that $g^{-1} F(g) \in \Pb_0\cdot
F(\Pb_0)$ 
if and only if $(gg_0)^{-1} F(gg_0) \in \Pb_I w \Pb_{F(I)}$). 

Let $\Qb$ be a parabolic subgroup of $\Gb$ containing $\Pb$. Let $J$ be its 
type. Then $I \subset J$, $\Qb=\lexp{g_0}{\Pb_J}$ and 
$g_0^{-1}F(g_0) \in \Pb_J w \Pb_{F(J)}$. Now, $\Qb$ is $F$-stable if and only 
if $F(J)=J$ and $w \in W_J$. This shows the equivalence of the two Theorems.

\bigskip

\noindent{\sc Remark 2 - } The condition {\it ``$W_Iw$ is not contained in 
a proper $F$-stable standard parabolic subgroup of $W$''} is equivalent to 
{\it ``$W_IwW_{F(I)}$ is not contained in 
a proper $F$-stable standard parabolic subgroup of $W$''}.

\bigskip

The rest of this paper is devoted to the proof of Theorem A'.
We fix a subset $I$ of $S$ and an element $w$ of $W$. We first 
recall two elementary facts. If $I \subset J$, let $\t_{IJ} : \PCB_I \to \PCB_J$ 
be the morphism of varieties that sends $\Pb \in \PCB_I$ to the unique 
parabolic subgroup of type $J$ containing $\Pb$. It is surjective. Moreover, 
\equat\label{tau ij}
\t_{IJ}(\Xb_I(w)) \subset \Xb_J(w)
\endequat
and
\equat\label{tau inverse}
\t_{IJ}^{-1}(\Xb_J(w)) = \bigcup_{W_I x W_{F(I)} \subset W_J w W_{F(J)}} \Xb_I(x).
\endequat

\bigskip

\noindent{\bfit First step: the ``only if'' part.} 
Assume that there exists a proper $F$-stable subset $J$ of $S$ such that 
$W_I w \subset W_J$. Then, by \ref{tau ij}, we have 
$\t_{IJ}(\Xb_I(w)) \subset \Xb_J(1)=\PCB_J^F$. 
Since $\Gb^F$ acts transitively on $\PCB_J^F$, we get 
$\t_{IJ}(\Xb_I(w))=\Xb_J(1)$. This shows that $\Xb_I(w)$ is not irreducible. 

\bigskip

\noindent{\bfit Second step: reduction to Borel subgroups.} 
By the previous step, we can concentrate on the ``if'' part. 
So, from now on, we assume that $W_I w$ is not contained 
in a proper $F$-stable parabolic subgroup of $W$. Then, by \ref{tau inverse}, 
we have 
$$\t_{\vide I}^{-1}(\Xb_I(w))=\bigcup_{x \in W_I w W_{F(I)}} \Xb(x).$$
Let $v$ denote the longest element of $W_I w W_{F(I)}$. Then 
every element $x$ of the double coset $W_I w W_{F(I)}$ satisfies 
$x \le v$ (here, $\le$ denotes the Bruhat order on $W$): this follows 
for instance from the fact that $\Pb_I w \Pb_{F(I)}$ is irreducible 
and is equal to $\cup_{x \in W_I w W_{F(I)}} \Bb w \Bb$. In particular, 
$v$ is not contained in a proper $F$-stable parabolic subgroup of $W$. 

Now, let $\Xb'=\bigcup_{x \in W_I w W_{F(I)}} \Xb(x)$. Then, 
since $\overline{\Xb(v)}=\bigcup_{x \le v} \Xb(x)$, we have 
$$\Xb(v) \subset \Xb' \subset \overline{\Xb(v)}.$$
So, since $\t_{\vide I}(\Xb')=\Xb_I(w)$, it is enough to show that $\Xb(v)$ 
is irreducible. 
In other words, we may, and we will, assume that $I=\vide$. 

\bigskip

\noindent{\bfit Third step: smooth compactification.} 
Let $(s_1,\dots,s_n)$ be a finite sequence of elements of $S$. Let 
\begin{multline*}
\Xbh(s_1,\dots,s_n)=\{(\Bb_1,\dots,\Bb_n) \in \BCB^n~|~
(\Bb_n,F(\Bb_1)) \in \overline{\OCB(s_n)} \\
\text{and } (\Bb_i,\Bb_{i+1}) \in \overline{\OCB(s_i)}
\text{ for }1\le i\le n-1\}.
\end{multline*}
If $\ell(s_1 \cdots s_n)=n$, then $\Xbh(s_1,\dots,s_n)$ is a smooth 
compactification of $\Xb(s_1\cdots s_n)$ (see \cite[Lemma 9.11]{delu}): 
in this case, 
\equat\label{irreductible}
\text{\it $\Xb(s_1\cdots s_n)$ is irreducible if and only 
if $\Xbh(s_1, \dots, s_n)$ is irreducible.}
\endequat

Note that $(\Bb,\dots,\Bb) \in \Xbh(s_1,\dots,s_n)$. We denote by 
$\Xbh^\circ(s_1,\dots,s_n)$ the connected (i.e. irreducible) component 
of $\Xbh(s_1,\dots,s_n)$ containing $(\Bb,\dots,\Bb)$. Let
$H(s_1,\dots,s_n)\subset \Gb^F$ 
be the stabilizer of $\Xbh^\circ(s_1,\dots,s_n)$. Let us now prove 
the following fact: 
\equat\label{stabilisateur}
\text{\it if $1 \le i_1 < \dots < i_r \le n$, then 
$H(s_{i_1},\dots,s_{i_r}) \subset H(s_1,\dots,s_n)$.}
\endequat

\noindent{\it Proof of \ref{stabilisateur} -} 
The map
$f : \Xbh(s_{i_1},\dots,s_{i_r}) \longto \Xbh(s_1,\dots,s_n)$
defined by 
\begin{multline*}
f(\Bb_1,\dots,\Bb_1) = \\
(\Bb_1,\dots,
\underbrace{\Bb_1}_{\substack{\text{$i_1$-th} \\ \text{position}}},
\Bb_2,\dots,
\underbrace{\Bb_{r-1}}_{\substack{\text{$i_{r-1}$-th} \\ \text{position}}},
\Bb_r,\dots,\underbrace{\Bb_r}_{\substack{\text{$i_r$-th} \\ \text{position}}},
F(\Bb_1),\dots,F(\Bb_1))
\end{multline*}
is a $\Gb^F$-equivariant morphism of varieties. Moreover, 
$$f(\underbrace{\Bb,\dots,\Bb}_{\text{$r$ times}})=
(\underbrace{\Bb,\dots,\Bb}_{\text{$n$ times}}).$$
In particular, $f(\Xbh^\circ(s_{i_1},\dots,s_{i_r}))$ is contained 
in $\Xbh^\circ(s_1,\dots,s_n)$. This proves the expected inclusion 
between stabilizers.\fin

\bigskip

\noindent{\bfit Last step: twisted Coxeter element.} 
The quotient variety
$$\Gb^F\backslash
\{g\in \Gb|g^{-1}F(g)\in \Bb w\Bb\}$$
is irreducible (it is isomorphic to $\Bb w \Bb$ through the Lang 
map $\Gb^F g \mapsto g^{-1} F(g)$), hence
$\Gb^F\backslash \Xb(w)$ is irreducible as well. So, 
\equat\label{3}
\text{\it $\Gb^F$ permutes transitively the irreducible components of $\Xb(w)$.}
\endequat

Let $w=s_1\cdots s_n$ be a reduced decomposition of $W$ as a product 
of elements of $S$.
By \ref{irreductible} and \ref{3}, it suffices to show that
$H(s_1,\dots,s_n)=\Gb^F$. 
Since $w$ does not belong to any $F$-stable proper 
parabolic subgroup of $W$, there exists a sequence $1 \le i_1 < \dots < i_r \le n$ 
such that $(s_{i_k})_{1 \le k \le r}$ is a family of representatives 
of $F$-orbits in $S$. By \ref{stabilisateur}, 
we have $H(s_{i_1},\dots,s_{i_r}) \subset H(s_1,\dots,s_n)$. 
But, by \cite[Proposition 4.8]{lucox}, $\Xb(s_{i_1}\cdots,s_{i_r})$ is
irreducible 
so, again by \ref{irreductible} and \ref{3}, $H(s_{i_1},\dots,s_{i_r})=\Gb^F$. 
Therefore, $H(s_1,\dots,s_n)=\Gb^F$, as expected. 

\bigskip

\noindent{\bf Acknowledgements.} We thank F. Digne and J. Michel for 
fruitful discussions on these questions. We thank P. Deligne for 
the clarification of the scope of validity of the Theorem.


\bigskip

\begin{thebibliography}{DIMIROU}

\bibitem[DeLu]{delu} P.~Deligne and G.~Lusztig, 
	{\em Representations of reductive groups over finite fields}, 
 	Ann. of Math. {\bf 103} (1976), 103--161.
\bibitem[DiMi2]{dimi} F.~Digne and J.~Michel, 
        {\em Endomorphisms of Deligne-Lusztig varieties}, 
        preprint (2005), {\tt math.RT/0509011}. 
\bibitem[Lu]{lucox} G.~Lusztig,
	{\em Coxeter orbits and eigenspaces of Frobenius},
	Inv. Math. {\bf 38} (1976), 101--159.
\end{thebibliography}
\end{document}